# NONEXISTENCE OF RANDOM GRADIENT GIBBS MEASURES IN CONTINUOUS INTERFACE MODELS IN $D = 2$

BY AERNOUT C. D. VAN ENTER AND CHRISTOF KÜLSKE

*University of Groningen*

We consider statistical mechanics models of continuous spins in a disordered environment. These models have a natural interpretation as effective interface models. It is well known that without disorder there are no interface Gibbs measures in infinite volume in dimension $d = 2$, while there are "gradient Gibbs measures" describing an infinite-volume distribution for the increments of the field, as was shown by Funaki and Spohn.

In the present paper we show that adding a disorder term prohibits the existence of such gradient Gibbs measures for general interaction potentials in $d = 2$. This nonexistence result generalizes the simple case of Gaussian fields where it follows from an explicit computation.

In $d = 3$ where random gradient Gibbs measures are expected to exist, our method provides a lower bound of the order of the inverse of the distance on the decay of correlations of Gibbs expectations w.r.t. the distribution of the random environment.

**1. Introduction.**

1.1. *The setup.* Our model is given in terms of the formal infinite-volume Hamiltonian

(1) $$H[\eta](\varphi) = \tfrac{1}{2} \sum_{i,j} p(i-j) V(\varphi_i - \varphi_j) - \sum_i \eta_i \varphi_i.$$

Here the fields (unbounded continuous spins) $\varphi_i \in \mathbb{R}$ represent height variables of a random surface at the site $i \in \mathbb{Z}^d$. Such a model is motivated as an effective model for the study of phase boundaries at a mesoscopic level in statistical mechanics. The disorder configuration $\eta = (\eta_i)_{i \in \mathbb{R}^d}$ denotes an arbitrary fixed configuration of external fields, modeling a "quenched" (or frozen) random environment.

For background and various earlier results about both continuous and discrete interface models without disorder, see [2, 5, 8, 12, 13] and references therein; for results about discrete interface models in the presence of disorder, see [3, 4].

*Assumptions.* The pair potential $V(t)$ is assumed to be even, $V(t) = V(-t)$, and continuously differentiable. We require that $V$ grows faster than linearly to infinity, that is, $\lim_{|t|\uparrow\infty} \frac{V(t)}{|t|^{1+\varepsilon}} = \infty$ for some positive $\varepsilon$. $p(\cdot)$ is the transition kernel of a simple random walk on $\mathbb{Z}^d$, assumed to be symmetric and of finite range. We further demand for simplicity that the random fields $\eta_i$ be i.i.d. under the distribution $\mathbb{P}$, symmetric and with finite nonzero second moment. We denote the expectation w.r.t. $\mathbb{P}$ by the symbol $\mathbb{E}$.

*Vector fields, fields and gradient fields.* We call the set of sites $\mathbb{Z}^d$ with oriented edges between $i, j$ whenever $p(i - j) > 0$ the *graph of the random walk*. We call a *vector field* (or vector field configuration) a map from the set of oriented edges $ij$ of this graph to $\mathbb{R}$ such that $V_{i,j} = -V_{j,i}$. Every field configuration $U = (U_i)_{i\in\mathbb{Z}^d}$ gives rise to a vector field configuration $U' = U'(U)$ by $U'_{i,j} = U(i) - U(j)$. In this case we call $U'$ *the gradient field of* $U$. There does not need to exist such a function for a general vector field $V$. Its existence is equivalent to the *loop condition* $V_{i,j} + V_{j,k} + V_{k,l} + V_{l,i} = 0$ along plaquettes $i, j, k, l$. We denote the set of all field configurations in infinite volume by $\Omega$ and the set of all gradient field configurations in infinite volume by $\Omega'$.

*Gibbs measures and gradient Gibbs measures.* The quenched finite-volume Gibbs measures (or local specification) corresponding to the Hamiltonian (1) in a finite volume $\Lambda \subset \mathbb{Z}^d$, a boundary condition $\hat\varphi$ and a fixed disorder configuration in $\Lambda$ are given by the standard expression

$$\int \mu_\Lambda^{\hat\varphi}[\eta](\varphi)(F(\varphi))$$

$$:= \int d\varphi_\Lambda F(\varphi_\Lambda, \hat\varphi_{\Lambda^c})$$

(2)
$$\times \exp\biggl\{-\tfrac{1}{2} \sum_{i,j\in\Lambda} p(i-j)V(\varphi_i - \varphi_j)$$

$$- \sum_{i\in\Lambda, j\in\Lambda^c} p(i-j)V(\varphi_i - \hat\varphi_j) + \sum_{i\in\Lambda} \eta_i\varphi_i\biggr\}$$

$$\times Z_\Lambda^{\hat\varphi}[\eta]^{-1},$$

where $Z_\Lambda^{\hat\varphi}[\eta]$ denotes the normalization constant that turns the last expression into a probability measure. It is a simple matter to see that the growth condition on $V$ guarantees the finiteness of the integrals appearing in (2) for all arbitrarily fixed choices of $\eta$.

We note that the Hamiltonian $H[\eta]$ changes only by a configuration-independent constant under the joint shift $\varphi_x \mapsto \varphi_x + c$ of all height variables with the same $c$. This holds true for any fixed configuration $\eta$. Hence, finite-volume Gibbs measures transform under a shift of the boundary condition by a shift of the integration variables. Using this invariance under height shifts we can lift the finite-volume measures to measures on gradient configurations, defining the *gradient finite-volume Gibbs measures* (gradient local specification). These are the probability kernels from $\Omega'$ to $\Omega'$ given by

$$(3) \qquad \int (\mu')_\Lambda^{\bar\varphi'}[\eta](d\varphi')G(\varphi') := \int \mu_\Lambda^{\bar\varphi}[\eta](d\varphi)G(\varphi'(\varphi)),$$

where $\bar\varphi$ is any field configuration whose gradient field is $\bar\varphi'$.

Finally, we call $\mu'$ a *gradient Gibbs measure* if it satisfies the DLR equation, that is,

$$(4) \qquad \int \mu'(d\bar\varphi') \int (\mu')_\Lambda^{\bar\varphi'}(d\varphi')f(\varphi') = \int \mu'(d\varphi')f(\varphi').$$

Here we dropped the $\eta$ from our notation.

Gradient Gibbs measures have the advantage that they may exist, even in situations where a proper Gibbs measure does not. When this happens, it means that the interface is locally smooth, although at large scales it has too large fluctuations to stay around a given height. For background about these objects, we refer to [8, 9, 12].

*Translation-covariant (gradient) Gibbs measures.* Denote by $\tau_r\varphi = (\varphi_{i-r})_{i\in\mathbb{Z}^d}$ the shift of field configurations on the lattice by a vector $r \in \mathbb{Z}^d$. A measurable map $\eta \mapsto \mu[\eta]$ is called a *translation-covariant random Gibbs measure* if $\mu[\eta]$ is a Gibbs measure for almost any $\eta$, and if it behaves appropriately under lattice shifts, that is, $\int \mu[\tau_r\eta](d\varphi)F(\varphi) = \int \mu[\eta](d\varphi)F(\tau_r\varphi)$ for all translation vectors $r$. This means that the functional dependence of the Gibbs measure on the underlying disordered environment is the same at every point in space. This notion goes back to [1] and generalizes the notion of a translation-invariant Gibbs measure to the set-up of disordered systems.

Naturally, a measurable map $\eta \mapsto \mu'[\eta]$ is called a *translation-covariant random gradient Gibbs measure* if $\int \mu'[\tau_r\eta](d\varphi')F(\varphi') = \int \mu'[\eta](d\varphi')F(\tau_r\varphi')$ for all $r$.

1.2. *Main results.* A main question to be asked in interface models is whether the fluctuations of an interface that is restricted to a finite volume will remain bounded when the volume tends to infinity, so that there is an infinite-volume Gibbs measure (or gradient Gibbs measure) describing a localized interface. This question is well understood in translation-invariant continuous-height models, and it is the purpose of this note to discuss such models in a random environment.

Let us start by discussing the nonrandom model in the physically interesting dimension $d=2$. We start with the Gaussian model where $V(t) = \frac{t^2}{2}$. Gaussian models are simple because explicit computations can be done. These show that the Gibbs measure $\mu(d\varphi)$ does not exist in infinite volume, but the gradient field (gradient Gibbs measure) does exist in infinite volume. Both statements may easily be derived by looking at the properties of (differences of) the matrix elements of $(I - P_\Lambda)^{-1}$ which appears here as a covariance matrix, where $P_\Lambda$ is the transition operator, restricted to $\Lambda$, given in terms of the random walk kernel $p$. Equivalently, one may say that the infinite-volume measure exists conditioned on the fact that one of the variables $\varphi_i$ is pinned at the value zero. Funaki and Spohn showed more generally that for convex potentials $V$ there are tilted gradient Gibbs measures and, conversely, a gradient Gibbs measure is uniquely determined by the tilt [8, 9, 12].

For (very) nonconvex $V$, new phenomena are appearing: There may be a first-order phase transition in the temperature where the structure of the interface (at zero tilt) changes, as shown by Biskup and Kotecky [2]. This phenomenon is related to the phase transition seen in rotator models with very nonlinear potentials exhibited in [6, 7]. The basic mechanism is an energy–entropy transition such as was first proven for the Potts model for a sufficiently large number of spin values [10].

What can we say for the random model? In [11] the authors showed a quenched deterministic lower bound on the fluctuations in the anharmonic model in a finite box of the order square root of the sidelength, uniformly in the disorder. In particular, this implies that there will not be any disordered infinite-volume Gibbs measures in $d=2$. This latter statement is not surprising since there is already nonexistence of the unpinned interface without disorder. But what will happen to the gradient Gibbs measure that is known to exist without disorder, once we allow for a disordered environment?

*Gaussian results and predictions.* Let us first look in some detail at the special case of a Gaussian gradient measure where $V(t) = \frac{t^2}{2}$, specializing to nearest neighbor interactions. Then, for any fixed configuration $\eta_\Lambda$, the finite-volume gradient Gibbs measure with zero boundary condition is a

Gaussian measure with expected value

$$X_{ij}^\Lambda[\eta] := \int \mu_\Lambda'[\eta](d\varphi')(\varphi_{ij}') = \sum_{y \in \Lambda} T_{ij,y}^\Lambda \eta_y$$

(5)

$$\text{where } T_{ij,y}^\Lambda = (-\Delta_\Lambda)_{i,y}^{-1} - (-\Delta_\Lambda)_{j,y}^{-1}.$$

The matrix $\text{cov}_\Lambda(\varphi_{ij}; \varphi_{lm}) = (-\Delta_\Lambda)_{i,l}^{-1} - (-\Delta_\Lambda)_{i,m}^{-1} - (-\Delta_\Lambda)_{j,l}^{-1} + (-\Delta_\Lambda)_{j,m}^{-1}$ is the same as in the model without disorder; its infinite-volume limit exists in any dimension $d \geq 2$, by a simple computation.

What about the infinite-volume limit of the mean value (5), as a function of the disorder? We note that $X^\Lambda[\eta] := (X_{ij}^\Lambda[\eta])_{\substack{ij \in \Lambda \times \Lambda \\ i \sim j}}$ is itself a random vector field with covariance

(6) $$C_\Lambda(ij, lm) = \mathbb{E}(\eta_0^2) \sum_{y \in \Lambda} T_{ij,y}^\Lambda T_{lm,y}^\Lambda.$$

Moreover, it is also a gradient vector field, since the loop condition carries over from $\varphi'$, by linearity of the Gibbs expectation.

For its variance, we have, in particular,

(7) $$C_\Lambda(ij, ij) = \mathbb{E}(\eta_0^2) \sum_{y \in \Lambda} (T_{ij,y}^\Lambda)^2.$$

In two dimensions the infinite-volume limit of this expression does not exist since $\int^N r(\frac{d}{dr} \log r)^2 \, dr \sim \log N$, when the sidelength $N$ of the box diverges to infinity. In dimension $d > 2$, we have $\int^N r^{d-1}(\frac{d}{dr} r^{-(d-2)})^2 \, dr \sim \int^N r^{-(d-1)} \, dr$, so the fluctuations stay bounded.

In particular, this explicit computation shows that in the Gaussian model there cannot be infinite-volume random gradient Gibbs measures in $d = 2$. Indeed, already the local—short-distance—fluctuations are roughening up the interface.

It is the main result of this paper to show that this result persists also for an anharmonic potential where explicit computations are not possible.

THEOREM 1.1 (Nonexistence in $d = 2$). *Suppose $d = 2$. Then there does not exist a translation-covariant random gradient Gibbs measure $\mu'[\eta](d\varphi')$ that satisfies the integrability condition*

(8) $$\mathbb{E}\left|\int \mu'[\eta](d\varphi') V'(\varphi_{ij}')\right| < \infty$$

*for sites $i \sim j$.*

Let us consider dimension $d = 3$ where translation-covariant infinite-volume gradient measures are believed to exist. The next result shows that, if they do, they must have slow decay of correlations w.r.t. the random environment.

THEOREM 1.2 (Slow decay of correlations in $d=3$). *Suppose that $d=3$ and that $\mu'[\eta](d\varphi')$ is a random gradient Gibbs measure. Put*

$$(9) \qquad C(ij,kl) := \mathbb{E}\Big(\int \mu'[\eta](d\varphi')V'(\varphi'_{ij}) \int \mu[\eta](d\varphi')V'(\varphi'_{kl})\Big)$$

*for sites $i \sim j$ and $k \sim l$.*
*Then*

$$(10) \qquad \lim_{r\uparrow\infty} \sup_{i,j,k,l\,:\,|i-k|\geq r} \frac{|C(ij,kl)|}{r^{-(1+\varepsilon)}} = \infty$$

*for all $\varepsilon > 0$.*

REMARK. Note that in $d=3$ for the case of the quadratic nearest neighbor potential $V(t) = \frac{t^2}{2}$ a translation-covariant random gradient Gibbs measure exists by explicit computation. It is the Gaussian field whose covariance and mean are given by the infinite-volume limits of (6) and (5). An explicit computation with this Gibbs measure also shows that the power of the bound given in the last theorem is optimal; see Section 2.3.

The method of proof of both theorems relies on a surface-volume comparison. It is inspired by the Aizenman–Wehr method (devised in [1] and used on discrete interfaces in [4]), but different in several aspects. As opposed to the Aizenman–Wehr Ansatz, where free energies are considered, we derive a discrete divergence equation $\eta = \nabla X$ where the external random field (the disorder) acts as a source and the vector field $X$ is provided by the expectation of $V'(\varphi_i - \varphi_j)$ with respect to the hypothesized gradient Gibbs measure. Exploiting a discrete version of Stokes' theorem and a volume versus surface comparison between terms, we arrive at a contradiction.

## 2. Proof of Theorems 1.1 and 1.2.

2.1. *Proof of Theorem* 1.1. The method of proof is to argue from the existence of a translation-covariant gradient Gibbs measure to a contradiction in $d=2$. To do this, we start with the following definition.

DEFINITION 2.1. Let $\mu'$ be an infinite-volume gradient Gibbs measure, in a random or nonrandom model. Then we call the vector field $X$ on the graph of the random walk given by

$$(11) \qquad X_{ij} := \int \mu'(d\varphi')(V'(\varphi'_{ij}))$$

the associated vector field.

This indeed defines a vector field because, by symmetry, $V'(x) = -V'(-x)$ and, hence, $X_{ij} = -X_{ji}$. Of course, the same definition can be made in finite volume, but for our proof, we will work immediately in infinite volume.

We note that a gradient Gibbs measure in the Gaussian model provides even a gradient vector field, since $V'(x) = cx$ is a linear function, and so the loop condition carries over from $\varphi'_{ij}$ to $X_{ij}$. For general nonquadratic potentials $V$, the vector fields $X_{ij}$ will not be gradient fields. This explains why, in the anharmonic model, we must work with vector fields (functions on the edges).

PROPOSITION 2.2. *Let $\mu'[\eta](d\varphi')$ be a random infinite-volume gradient Gibbs measure and $X_{ij}[\eta] := \int \mu'[\eta](d\varphi')(V'(\varphi'_{ij}))$ the associated vector field on the random walk graph for $i \sim j$. Then*

(12) $$\eta_i = \sum_j p(j-i) X_{ij}[\eta]$$

*for all $i \in \mathbb{Z}^d$. This equation will be called the "divergence equation".*

REMARK. This equation is a discrete version of the equation $\eta = \nabla \cdot X$ which holds in continuous space $\mathbb{R}^d$.

Comparing with Maxwell's equations, $X$ plays the role of the electric field, and $\eta$ plays the role of the electric charge. In the Gaussian case, we know that $X$ is a gradient-field, that is, curl-free. In the case of general potentials $V$ it will not be. Note that the kernel of this equation, that is, the vector fields $X$ which are the solutions of the equation $0 = \nabla \cdot X$, consists of all divergence-free vector fields, which are plenty. In particular, the divergence equation does not allow to determine $X$ uniquely in terms of $\eta$.

PROOF OF PROPOSITION 2.2. We look at the one-site local specification at the site $i$. Take a single-site integral over $\varphi_i$ appearing in the partition function and use partial integration to write

(13)
$$\int d\varphi_i \exp\left(\sum_{j \sim i} p(j-i) V(\varphi_i - \varphi_j)\right)(\eta_i \exp(\eta_i \varphi_i))$$
$$= \int d\varphi_i \left(-\frac{\partial}{\partial \varphi_i} \exp\left(-\sum_{j \sim i} p(j-i) V(\varphi_i - \varphi_j)\right)\right) \exp(\eta_i \varphi_i)$$
$$= \int d\varphi_i \sum_{j \sim i} p(j-i) V'(\varphi_i - \varphi_j)$$
$$\times \exp\left(-\sum_{j \sim i} p(j-i) V(\varphi_i - \varphi_j) + \eta_i \varphi_i\right).$$

Now, integrate over the remaining $\varphi_k$ for $k$ in a finite volume $\Lambda$ to see that

$$\eta_i = \sum_j p(j-i) \int \mu_\Lambda^{\bar\varphi}[\eta](d\varphi_\Lambda) V'(\varphi_i - \varphi_j)$$

(14)

$$= \sum_j p(j-i) \int \mu_\Lambda^{\bar\varphi'}[\eta](d\varphi'_\Lambda) V'(\varphi'_{ij}).$$

Integrating this equation over the boundary condition $\bar\varphi'$ w.r.t. $\mu'[\eta]$ and using the DLR equation for the gradient measure implies the proposition. □

Summing the divergence equation over $i$ in a finite volume $\Lambda$, we note that the contributions of the edges that are contained in $\Lambda$ vanish, due to the property of $X$ being a vector field. Hence, we arrive at the following corollary.

COROLLARY 2.3 (Integral form of divergence equation).

(15) $$\sum_{i \in \Lambda} \eta_i = \sum_{\substack{i \sim j \\ i \in \Lambda, j \in \Lambda^c}} p(j-i) X_{ij}[\eta].$$

Note that the sum of boundary terms on the r.h.s. plays the role of a surface integral in the Stokes equation.

Let us now specialize to $d = 2$ and prove Theorem 1.1. We choose $\Lambda = \{-L, -L+1, \ldots, L\}^2$ and normalize by $\frac{1}{L}$.

We remark that the sum over boundary bonds $ij$ decomposes into the 4 sides of a square, whose bonds will be denoted by $B_L(\nu)$, $\nu = 1, 2, 3, 4$.

So we have

$$\frac{1}{L} \sum_{i \in \Lambda} \eta_i = \frac{1}{L} \sum_{i \in \Lambda, j \in \Lambda^c} p(i-j) X_{ij}[\eta]$$

(16)

$$= \sum_{\nu=1,2,3,4} \frac{1}{L} \sum_{\langle i,j \rangle \in B_L(\nu)} p(i-j) X_{ij}[\eta].$$

By the ergodic theorem, each of the four sums converges to its expected value

(17) $$\lim_{L \uparrow \infty} \frac{1}{L} \sum_{\langle i,j \rangle \in B_L(\nu)} p(i-j) X_{ij}[\eta] = 2 \sum_{j : j \cdot e_\nu > 0} p(j) \mathbb{E}(X_{0j}[\eta])$$

in whatever sense (e.g., almost surely or in $L^2$), where $e_1, e_2, e_3 = -e_1, e_4 = -e_2$ are unit vectors on $\mathbb{Z}^d$.

By the CLT, the l.h.s. of equation (16) does not converge almost surely to a constant; indeed, $\frac{1}{L}\sum_{i\in\Lambda}\eta_i$ converges only in distribution to a Gauss distribution with a nonzero variance. This is a contradiction to convergence to a delta-distribution, which proves Theorem 1.1.

2.2. *Proof of Theorem 1.2.* Take the square of the integral form of the divergence equation and take its expectation w.r.t. the measure $\mathbb{P}$:

$$\text{(18)} \qquad \mathbb{E}(\eta_0^2)|\Lambda| = \sum_{\substack{i\sim j \\ i\in\Lambda, j\in\Lambda^c}} \sum_{\substack{k\sim l \\ k\in\Lambda, l\in\Lambda^c}} p(j-i)p(k-l)C_{ij,kl}.$$

Assume the bound

$$\text{(19)} \qquad C_{ij,kl} \leq \text{Const}(1+|i-k|^2)^{-q}.$$

Let us take for $\Lambda$ a ball w.r.t. the Euclidean metric of radius $L$. Then, for large $L$, the l.h.s. behaves like a constant times $L^3$.

The large-$L$ asymptotics of the r.h.s. of (18) is provided by the large-$L$ behavior of the double integral

$$\text{(20)} \qquad \int_{LS^2} d\lambda(x) \int_{LS^2} d\lambda(y)(1+|x-y|^2)^{-q},$$

where $LS^2$ denotes the two-dimensional sphere with radius $L$ and $\lambda$ is the Lebesgue measure on the sphere. By rotation-invariance, this equals

$$\text{(21)} \quad cL^2 \int_{LS^2} d\lambda(y)(1+|y-Le|^2)^{-q} = cL^4 \int_{S^2} d\lambda(z)(1+L^2|z-e|^2)^{-q},$$

where $e$ is the North Pole of $S^2$. Using polar coordinates $\cos\theta = s$ for the point $z$, where $s=1$ is corresponding to the point $e$, we have $|z-e|^2 = 2(1-s)$. This gives

$$\text{(22)} \quad \begin{aligned} &\int_{S^2} d\lambda(z)(1+L^2|z-e|^2)^{-q} \\ &= 2\pi \int_{-1}^{1} ds(1+L^2 2(1-s))^{-q} = \frac{2\pi 2^{1-q}}{1-q}(1+2L^2)^{-q}. \end{aligned}$$

By (18), we thus have under the assumption of (19) that $L^3 \leq cL^{4-2q}$, for large $L$, which can only be true for $q \leq \frac{1}{2}$. This implies that $C_{ij,kl}$ can not decay faster than the inverse distance between $i$ and $k$.

2.3. *Sharpness of polynomial decay in $d=3$.*

PROPOSITION 2.4. *Suppose $d=3$. Let $\mu'[\eta](d\varphi')$ be a random gradient Gibbs measure for the Gaussian model $V(x) = \frac{x^2}{2}$ with i.i.d. $\eta_i$ with finite second moment.*

*Then*

$$\text{(23)} \qquad \lim_{r\uparrow\infty} \sup_{i,j,k,l:|i-k|\geq r} \frac{|C(ij,kl)|}{r^{-1}} < \infty.$$

REMARK. This shows that the fluctuation lower bound is sharp in $d=3$.

PROOF OF PROPOSITION 2.4. Let us choose $i = -Re$ and $k = Re$, where $e$ is a unit coordinate vector so $2R$ is the distance between $i$ and $k$. To estimate the large-$R$ asymptotics of the decay of the covariance in infinite volume which is given by

$$\text{(24)} \qquad C_{\mathbb{Z}^3}(ij,lm) = E(\eta_0^2) \sum_{y\in\Lambda} T^{\mathbb{Z}^3}_{ij,y} T^{\mathbb{Z}^3}_{lm,y},$$

we are led to consider the large-$R$ asymptotics of the following integral:

$$\text{(25)} \qquad I(R) := \int d^3 y \frac{1}{1+|y-Re|^2} \frac{1}{1+|y+Re|^2}.$$

Now the proposition follows by an explicit computation.

Indeed, using cylindrical coordinates, this integral can be rewritten as

$$
\begin{aligned}
I(R) &= 2\pi \int_0^\infty dz \int_0^\infty d\bar{r}\,\bar{r}\, \frac{1}{1+(z-R)^2+\bar{r}^2} \frac{1}{1+(z+R)^2+\bar{r}^2} \\
\text{(26)} \quad &= \frac{\pi}{4R} \int_0^\infty \frac{dz}{z} \lim_{S\uparrow\infty} \int_0^S ds \left( \frac{1}{1+(z-R)^2+s} - \frac{1}{1+(z+R)^2+s} \right) \\
&= \frac{\pi}{4R} J(R).
\end{aligned}
$$

Here we have used the substitution $s = r^2$ and

$$\text{(27)} \qquad J(R) := \int_0^\infty \frac{dz}{z} \log \frac{1+(z+R)^2}{1+(z-R)^2} = \int_0^\infty \frac{du}{u} \log \frac{R^{-2}+(u+1)^2}{R^{-2}+(u-1)^2}.$$

Finally, the function $J(R)$ is increasing in $R$ and has the finite limit

$$\text{(28)} \qquad \lim_{R\uparrow\infty} J(R) = \int_0^\infty \frac{du}{u} \log \frac{(u+1)^2}{(u-1)^2} = 8 \int_0^1 \frac{\log x}{x^2-1} dx = \pi^2.$$

(See an integral table or check with Mathematica.) This shows the proposition. □

DEPARTMENT OF MATHEMATICS
AND COMPUTING SCIENCES
UNIVERSITY OF GRONINGEN
BLAUWBORGJE 3
9747 AC GRONINGEN
THE NETHERLANDS
E-MAIL: aenter@phys.rug.nl
        kuelske@math.rug.nl
URL: http://www.statmeca.fmns.rug.nl/
     http://www.math.rug.nl/~kuelske/